\documentclass[a4paper, 10pt, oneside, onecolumn, preprint]{elsarticle}
\usepackage[english]{babel}
\usepackage{a4wide,color,graphicx,amsmath,amssymb,verbatim,mathrsfs,latexsym}
\pdfoutput=1
\allowdisplaybreaks
\newtheorem{theo}{Theorem}
\newtheorem{prop}{Proposition}

\newtheorem{claim}{Claim}

\newcommand{\R}{\mathbb{R}}

\newcommand{\N}{\mathbb{N}}

\newcommand{\ep}{\varepsilon}
\newcommand{\pa}{\partial}
\newcommand{\Div}{\textrm{div}\,}
\newcommand{\Curl}{\textrm{curl}\,}
\newcommand{\fun}[5]{  
\begin{array}{cccc}
{#1}\,: & #2 & \to & #3 \\
\phantom{{#1}\,:\,} & #4 & \mapsto & #5
\end{array}
}
\newcommand{\bn}{\overline{n}}
\newcommand{\bV}{\overline{V}}
\newcommand{\OmT}{\Omega\times (0,T)}
\title{Analysis of a drift-diffusion model with velocity saturation for spin-polarized transport in semiconductors}
\author{Nicola Zamponi}
\address{Institute for Analysis and Scientific Computing, Vienna University of  
	Technology, Wiedner Hauptstra\ss e 8--10, 1040 Wien, Austria}
\ead{nicola.zamponi@tuwien.ac.at}
\date{\today}

\begin{document}

\begin{frontmatter}
 
\begin{abstract}
A system of drift-diffusion equations with electric field under Dirichlet boundary conditions is analyzed. 
The system of strongly coupled parabolic equations for particle density and spin density vector describes 
the spin-polarized semi-classical electron transport in ferromagnetic semiconductors.
The presence of a nonconstant and nonsmooth magnetization vector, solution of the Landau-Lifshitz equation,
causes the diffusion matrix to be dependent from space and time and to have in general poor regularity properties, thus making the analysis challenging.
To partially overcome the analytical difficulties the velocity saturation hypothesis is made, which results in a bounded drift velocity.
The global-in-time existence and uniqueness of weak solutions is shown by means of a semi-discretization in time, which
yields an elliptic semilinear problem, and a quadratic entropy inequality, which allow for the limit of vanishing time step size. 
The convergence of the weak solutions to the steady state, under some restrictions on the parameters and data, is shown. 
Finally the higher regularity of solutions for a smooth magnetization in two space dimensions is shown through a diagonalization argument,
which allows to get rid of the cross diffusion terms in the fluid equations, and the iterative application of 
Gagliardo-Nirenberg inequalities and a generalized version of Aubin lemma.
\end{abstract}


%

\end{frontmatter}


\section{Introduction}\label{sec.intro}

In this paper, we study global existence, uniqueness, long-time behavior and regularity of the solutions to the following 
drift-diffusion equations:\footnote{We adopt the Einstein convention of sum over repeated indexes.}
\begin{align}
 & \pa_t n_0 - \pa_i\left(\frac{D}{\eta^2}(J_0^i-p m_s J_s^i)\right) = 0\qquad\textrm{in }\OmT,\label{eq.n0}\\
 & \pa_t n_k - \pa_i\left(\frac{D}{\eta^2}\left(- p m_k J_0^i + (\eta\delta_{ks} + (1-\eta)m_k m_s) J_s^i \right)\right)\nonumber\\
 & \qquad -2\gamma\ep_{ijk} n_i m_j + \frac{n_k}{\tau} = 0\qquad (k=1,2,3)\qquad\textrm{in }\OmT,\label{eq.nk}\\
 & J_0^i = \pa_i n_0 - v_i n_0,\quad J_s^i = \pa_i n_s - v_i n_s ,\qquad(1\leq i,s\leq 3)\quad\textrm{in }\OmT,\label{eq.J}
\end{align}
where $\pa_t = \pa/\pa t$, $\pa_i = \pa/\pa x_i$ are the partial derivatives with respect to $t$ and $x_i$, for $i=1,2,3$, respectively,
$\Omega\subset\R^d$ with $d\leq 3$ is a bounded domain, $D(x)$ is the diffusion coefficient, $p(x)\in (-1,1)$ is the spin polarization, 
$\eta(x) = \sqrt{1-p(x)^2}$, $\vec v(x,t) = (v_1(x,t),v_2(x,t),v_3(x,t))$ is the scaled drift velocity,
$\vec m(x,t) = (m_1(x,t),m_2(x,t),m_3(x,t))$ is the magnetization vector and $\ep_{ijk}$ is the Levi-Civita tensor.
The system describes the evolution of a spin-polarized electron charge distribution with charge density $n_0(x,t)$, 
spin density vector $\vec n(x,t) = (n_1(x,t),n_2(x,t),n_3(x,t))$ in a ferromagnetic semiconductor.
The drift velocity is related to the electric field $-\nabla V(x,t)$ by the following relation \cite{Marder,YuCardona}:
\begin{align}
& v_i = -\mu(|\nabla V|)\pa_i V \qquad (i=1,2,3), \label{v} 
\end{align}
where the function $\mu\,:\, [0,\infty)\to (0,\infty)$ is the so-called {\em electron mobility}, and 
the electric potential $V(x,t)$ is self-consistently given by the Poisson equation \cite{Jue09}:
\begin{align}
 & -\lambda_D^2\Delta V = 2 n_0 - C\qquad\textrm{in }\OmT,\label{eq.poi}
\end{align}
where $C = C(x)$ is the doping profile and $\lambda_D>0$ is the (constant) scaled Debye length. \\
For physical reasons, we require that $\mu$ satisfies:
\begin{align}
 & \mu\in C^1([0,\infty)),\quad \mu'(s)\leq 0,\quad s\mu(s)\leq v^{sat},
 \quad |s\mu(s)-\overline{s}\mu(\overline{s})|\leq L|s-\overline{s}|
 \quad \forall s, \overline{s}\geq 0, \label{hp.mu}
\end{align}
for some constants $v^{sat}>0$ (called {\em saturation velocity}), $L>0$. 
Eqs.~\eqref{v}, \eqref{hp.mu} imply that the scaled drift velocity $v$ of the electrons, which is parallel to the 
electric field $-\nabla V$, remains bounded even when $|\nabla V|$ is large. 
This fact is in agreement with the experimental evidence that the charge carriers velocity in a material cannot be higher
than the saturation velocity due to collisions with the acoustic phonons \cite{YuCardona}. 
We impose the following Dirichlet boundary conditions on $n_0$, $\vec n$, $V$, and initial conditions on $n_0$, $\vec n$:
\begin{align}
 & n_0 = n_0^D,\quad \vec n = \vec n^D,\quad V = V^D,\qquad\textrm{on }\pa\OmT,\label{bc}\\
 & n_0(\cdot,0) = n_0^I,\quad \vec n(\cdot,0) = \vec n^I\qquad\textrm{on }\Omega.\label{ic}
\end{align}
The magnetization vector $\vec m$ has constant modulus and satisfies the Landau-Lifshitz equation \cite{Aharoni, Chikazumi}:
\begin{align}
 &\pa_t\vec m = \vec m\wedge\Delta\vec m - \vec m\wedge (\vec m\wedge\Delta\vec m),
 \qquad |\vec m|\equiv 1 \qquad\textrm{in }\OmT. \label{LLG.eq}
\end{align}
We impose the following homogeneous Neumann boundary conditions and initial conditions on $\vec m$ \cite{AlSoy91}:
\begin{align}
& \pa_\nu\vec m = 0 \qquad\textrm{on }\pa\OmT,\qquad
\vec m(\cdot,0) = \vec m^{I}\qquad \textrm{on }\Omega.\label{bic.m}
\end{align}
The diffusion matrix $A = (a_{ij})_{i,j=0,1,2,3}$ of the system,
\begin{equation}
 A = D\eta^{-2}\begin{pmatrix}
1 &  -p\vec m^T \\
-p\vec m & \eta I + (1-\eta) \vec m\otimes\vec m
\end{pmatrix} ,
\end{equation}
which depends on $(x,t)$ through $\vec m$ and $p$, is symmetric and positive definite. In fact, its spectral decomposition is:
\begin{equation}
 A = \frac{D}{1+p}\Pi^+ + \frac{D}{1-p}\Pi^- + \frac{D}{\eta}\Pi^\perp , \label{A.dec}
\end{equation}
where:
\begin{equation}
 \Pi^\pm = \frac{1}{2}
 \begin{pmatrix}
  1 & \pm\vec m^T \\
  \pm\vec m & \vec m\otimes\vec m
 \end{pmatrix},\qquad
 \Pi^\perp = 
 \begin{pmatrix}
  0 & 0\\
  0 & I-\vec m\otimes\vec m
 \end{pmatrix},
\end{equation}
are the projections operators in the eigenspaces related to the eigenvalues of $A$.
In particular, the eigenvalues of $A$, namely $D/(1+p)$, $D/(1-p)$, $D/\eta$, do not depend on $\vec m$ and are
thus time-independent.

\section{Theoretical considerations}\label{sec.theory}

The derivation of system \eqref{eq.n0}--\eqref{eq.J} from a linear Boltzmann equation with matrix collision operator can be found in \cite{NegPos}. 
A similar model has been considered in \cite{Jou} (Chapter 6.2): the system of parabolic PDEs for 
particle density and spin vector given by eq.~(6.2.2) is formally identical to eqs.~\eqref{eq.n0}--\eqref{eq.J} with $\eta=1$. Moreover, in the 
model presented in \cite{Jou}, the magnetization is given by a Landau-Lifshitz equation with an additional term, proportional to the cross product
between the spin and the magnetization, which provides a weak coupling between the Landau-Lifshitz equation and the drift-diffusion equations.
Other drift-diffusion models for spin-polarized electron transport have been derived in \cite{ElHajj} 
from a spinor Boltzmann equation with spin-flip and non spin-flip collision operators
under the hypothesis that the former is a small perturbation of the latter, and in \cite{Sai04} from a Wigner equation 
with semiclassical matrix collision operator under the assumption that spin-orbit coupling is small in comparison with the electron kinetic energy.
Both models do not contain any cross diffusion term and consist in a set of decoupled linear parabolic equations for the particle density and the 
spin density vector. 

Eqs.~\eqref{eq.n0}--\eqref{eq.J} have been analytically studied in \cite{JueNegShp} in the case of 
given constant magnetization $\vec m$ and mobility $\mu\equiv 1$. In such a case the cross diffusion terms in eqs.~\eqref{eq.n0}--\eqref{eq.J} 
can be removed, thus significantly weakening the coupling between equations, by considering the variables:
\begin{align}
 & n_+ = n_0 + \vec m\cdot\vec n , \qquad
 n_- = n_0 - \vec m\cdot\vec n , \qquad
 \vec n_\perp = \vec n - (\vec n\cdot\vec m)\vec m .\label{newvar}
\end{align}
The quantities $n_+$ and $n_-$ are referred to in \cite{JueNegShp} as spin-up and spin-down densities, respectively.
In particular, the equations for $n_+$, $n_-$ and $\vec n_\perp$ are decoupled:
\begin{align}
 & \pa_t n_\pm - \pa_i\left(\frac{D}{1\pm p}(\pa_i n_\pm + n_\pm\pa_i V)\right) \pm \frac{n_+ - n_-}{2\tau} = 0, \label{JNS.npm}\\
 & \pa_t \vec n_\perp - \pa_i\left(\frac{D}{\eta}(\pa_i \vec n_\perp + \vec n_\perp\pa_i V)\right) - 2\gamma\vec n_\perp\wedge\vec m
 + \frac{\vec n_\perp}{\tau} = 0.\label{JNS.nperp}
\end{align}
Thanks to the simple structure of eqs.~\eqref{JNS.npm}, \eqref{JNS.nperp}, 
the existence of bounded solutions to eqs.~\eqref{eq.n0}--\eqref{eq.poi}, with positive densities $n_\pm$, is shown in \cite{JueNegShp}
under the assumption that the diffusion coefficient $D$ and the polarization $p$ are constant. The same result is proved for system 
\eqref{eq.n0}--\eqref{v} when $D$ and $p$ are bounded functions of $x$, and $V = V(x)$ is a given function with bounded gradient.

The mathematical analysis of eqs.~\eqref{eq.n0}--\eqref{bic.m} is challenging because of the presence of cross diffusion terms 
with nonconstant and nonsmooth coefficients, unlike the problem considered in \cite{JueNegShp}. 
Since in this paper the magnetization $\vec m = \vec m(x,t)$ is given by eq.~\eqref{LLG.eq}, we cannot rely on the above argument to prove 
existence, boundedness and positivity of the solutions to eqs.~\eqref{eq.n0}--\eqref{ic}. 
In fact, for a nonconstant magnetization the right-hand sides of eqs.~\eqref{JNS.npm}, \eqref{JNS.nperp} are nonzero and the equations are coupled.
As a consequence, results of positivity or boundedness for $n_\pm$ cannot be achieved with usual methods (like e.g. Stampacchia truncation technique 
\cite{ZamJue13}).
Moreover, since in the case at hand eqs.~\eqref{JNS.npm}, \eqref{JNS.nperp} would depend also on the derivatives of $\vec m$, they could not be exploited 
to prove existence of solutions for system \eqref{eq.n0}--\eqref{ic} unless strong regularity assumptions on $\vec m$ are assumed.
Also a constant mobility, which has been treated in \cite{JueNegShp}, creates huge problems to the mathematical study of eqs.~\eqref{eq.n0}--\eqref{ic} 
in presence of nonconstant cross diffusion terms, because the particular (quadratic) structure of the drift terms $n_j\nabla V$ ($0\leq j\leq 3$), with 
$V$ given by the Poisson equation \eqref{eq.poi}, cannot be exploited with standard techniques to derive entropy estimates 
(see e.g. \cite{GaGr86,JueNegShp,ZamJueCMS}). 
The (physically justified) velocity saturation assumption allows us to overcome this difficulty, 
by providing an upper bound for the drift velocity modulus. 
We defer the study of eqs.~\eqref{eq.n0}--\eqref{eq.poi} with a constant mobility to a future work.

\section{Main results}\label{sec.mainresults}

Now we will state our main results and explain the ideas of the proofs.\\
Let $\pa\Omega\in C^1$. To fix a convenient notation let us define: 
$$n = (n_0,n_1,n_2,n_3),\quad n^D = (n_0^D,n_1^D,n_2^D,n_3^D),\quad n^I = (n_0^I,n_1^I,n_2^I,n_3^I), $$
$$ B = (b_{ij})_{i,j=0,1,2,3} , \qquad b_{ij} = 
\begin{cases} 
2\gamma\ep_{ijk}m_k - \tau^{-1}\delta_{ij} & 1\leq i\leq 3,\, 1\leq j\leq 3,\\
0 & \mbox{otherwise}.
\end{cases}$$
We can then rewrite eqs.~\eqref{eq.n0}--\eqref{eq.poi}, \eqref{bc}, \eqref{ic} in the following synthetic form:
\begin{align}
 &\pa_t n_i = \Div (a_{ij}(\nabla n_j - n_j v)) + b_{ij}n_j\qquad (0\leq i\leq 3)\quad\textrm{in }\OmT,\label{eq.n}\\
 & v = -\mu(|\nabla V|)\nabla V\quad\textrm{in }\OmT,\label{eq.v}\\
 & -\lambda_D^2\Delta V = 2 n_0 - C\qquad\textrm{in }\OmT,\label{eq.poi.bis}\\
 & n = n^D,\quad V=V^D\qquad\textrm{on }\pa\OmT,\label{eq.bc}\\
 & n(\cdot,0) = n^I\qquad\textrm{on }\Omega.\label{eq.ic}
\end{align}
For the analysis of eqs.~\eqref{eq.n}--\eqref{eq.ic} we exploited the following result concerning the solutions of the Landau-Lifshitz equation, 
which is the content of Theorem 1.5 in \cite{AlSoy91} and Theorem 1.4 in \cite{CarbouFabrie}.
\begin{prop}\label{prop.LLG}
For all $\vec m^I\in H^1(\Omega)$ with $|\vec m^I|\equiv 1$ in $\Omega$, a solution $\vec m$ to eqs.~\eqref{LLG.eq}, \eqref{bic.m} exists such that:
\begin{equation}
\vec m\in L^\infty(0,T; H^1(\Omega))\cap H^1(0,T; L^2(\Omega))\qquad\forall T>0.\label{reg.m} 
\end{equation}
Moreover if $\Omega\subset\R^2$, then a number $r>0$ exists such that, for all
$\vec m^I\in H^2(\Omega)$ with $|\vec m^I|\equiv 1$ in $\Omega$, $\pa_\nu\vec m^I = 0$ on $\pa\Omega$, $\|\nabla\vec m^I\|_{H^1(\Omega)}\leq r$,
a unique solution to eqs.~\eqref{LLG.eq}, \eqref{bic.m} exists satisfying:
\begin{equation}
 \vec m\in C^0([0,T]; H^2(\Omega))\cap L^2(0,T; H^3(\Omega))\cap C^1([0,T]; L^2(\Omega))\qquad \forall T>0 . \label{reg.m.2}
\end{equation}
\end{prop}
The first result we present is the following existence and uniqueness theorem.
\begin{theo}[Existence and uniqueness]\label{theo.ex} 
Let $\lambda_D>0$, $\gamma>0$, $T>0$. Moreover let us assume that $\vec m$ satisfies eq.~\eqref{reg.m} and
\begin{align}
 & n^{I}_j\in L^2(\Omega),\quad n^D_j,\, V^D\in H^1(\Omega),\quad C\in L^\infty(\Omega)\qquad (0\leq j\leq 3),\label{ex.hp.1}\\
 & D\in L^\infty(\Omega),\quad \inf_\Omega D>0,\quad \sup_\Omega |p| < 1.\label{ex.hp.2}
\end{align}
Then problem \eqref{eq.n}--\eqref{eq.ic} has a unique solution $(n,V)$ satisfying:
\begin{align}
 & n_j\in L^2(0,T; H^1(\Omega))\cap H^1(0,T; H^{-1}(\Omega))\cap L^\infty(0,T; L^2(\Omega))\qquad (0\leq j\leq 3),\label{reg.n}\\
 & V\in L^2(0,T; H^2(\Omega))\cap H^1(0,T; L^2(\Omega))\cap L^\infty(0,T; H^1(\Omega)).\label{reg.V}
 \end{align}
\end{theo}
The idea for the proof is to discretize eqs.~\eqref{eq.n}--\eqref{eq.ic} in time with the implicit Euler method with time step $h>0$. 
The problem we obtain is (here we neglect the boundary conditions \eqref{eq.bc} for the sake of simplicity):
\begin{align}
 & n_i^{(k)} - h\Div (a_{ij}^{(k)}\nabla n_j^{(k)}) 
  = n_i^{(k-1)} - h\Div (a_{ij}^{(k)}n_j^{(k)}v^{(k)})
 + h b_{ij}^{(k)}n_j^{(k)}\quad (0\leq i\leq 3)\quad\textrm{in }\Omega,\nonumber\\
 & v^{(k)} = -\mu(|\nabla V^{(k)}|)\nabla V^{(k)}\quad\textrm{in }\Omega,\quad
 -\lambda_D^2\Delta V^{(k)} = 2 n_0^{(k)} - C\qquad\textrm{in }\Omega,\nonumber
\end{align}
which can be solved by a standard fixed point argument. To take the limit $h\to 0$, thus solving the original problem \eqref{eq.n}--\eqref{eq.ic},
we derive a discrete entropy inequality for the following discrete quadratic entropy functional:
$$ S[n^{(k)}] = \frac{1}{2}\int_\Omega |n^{(k)}-n^D|^2 . $$
The inequality has the form:
\begin{align}
 & h^{-1}(S[n^{(k)}]-S[n^{(k-1)}]) + c_0\int_\Omega |\nabla n^{(k)}|^2 \leq c(S[n^{(k)}]+1),\nonumber
\end{align}
and yields gradient estimates for the solution, which allow us to take the limit for vanishing step size and so to obtain a solution of 
eqs.~\eqref{eq.n}--\eqref{eq.ic}. 
The uniqueness of solutions is achieved by deriving an entropy inequality for the relative entropy 
$$ S[n,\bn] = \frac{1}{2}\int_\Omega |n-\bn|^2 , $$
where $n$, $\bn$ are two solutions with the same initial and boundary data. Such inequality reads:
$$
\frac{d}{dt}S[n,\bn] \leq c(1+\|\bn\|_{H^1(\Omega)}^2)S[n,\bn]\qquad t>0 .
$$
Since $S[n,\bn](t=0)=0$, $S[n,\bn]\in L^\infty(0,T)$ and $\|\bn\|_{H^1(\Omega)}^2\in L^1(0,T)$, 
the above estimate allows to deduce $n=\bn$ once that the Gronwall lemma is applied.

The second result we present concerns the behaviour for $t\to\infty$ of the solution to eqs.~\eqref{eq.n}--\eqref{eq.ic}.
The steady state solution $(n^{eq}, V^{eq})$ for the system is defined by:
\begin{align}
 & \nabla n_0^{eq} + n_0^{eq}\mu(|\nabla V^{eq}|)\nabla V^{eq} = 0\qquad\textrm{in }\Omega,\label{lt.hp}\\
 & \lambda_D^2\Delta V^{eq} + 2 n_0^{eq} - C = 0 \qquad\textrm{in }\Omega,\label{lt.poi}\\
 & n_0^{eq} = n_0^D,\quad V^{eq} = V^D \qquad\textrm{on }\pa\Omega ,\label{lt.hp.2}\\
 & \vec n^{eq} = 0 \qquad\textrm{in }\Omega .\label{lt.hp.3}
\end{align}
We point out that, if $u^{eq}\equiv \log n_0^{eq}$, eqs.~\eqref{lt.hp}--\eqref{lt.hp.2} can be (formally) rewritten as:
\begin{align}
 & \Delta u^{eq} + \Div( \mu(|\nabla V^{eq}|)\nabla V^{eq} ) = 0\qquad\textrm{in }\Omega,\label{lt.pb.1}\\
 & \lambda_D^2\Delta V^{eq} + 2\exp(u^{eq}) - C = 0\qquad\textrm{in }\Omega,\label{lt.pb.poi}\\
 & u^{eq} = u^D\equiv\log n_0^D,\quad V^{eq} = V^D \qquad\textrm{on }\pa\Omega ,\label{lt.pb.bc}\\
 & \Curl(\mu(|\nabla V^{eq}|)\nabla V^{eq}) = 0\qquad\textrm{in }\Omega.\label{lt.pb.2}
\end{align}
While eqs.~\eqref{lt.pb.1}--\eqref{lt.pb.bc} can be solved with standard techniques for nonlinear elliptic equations, eq.~\eqref{lt.pb.2}
constitutes a nonlinear constraint for the solution $(n^{eq},V^{eq})$ of eqs.~\eqref{lt.pb.1}--\eqref{lt.pb.bc}, thus making the analytical study of 
the steady state problem tricky. For this reason, in this paper we just assume that there is a steady state $(n^{eq},V^{eq})$ and 
prove the convergence of the solution $(n,V)$ of eqs.~\eqref{eq.n}--\eqref{eq.ic} to $(n^{eq},V^{eq})$ as $t\to\infty$.
We defer the study of eqs.~\eqref{lt.pb.1}--\eqref{lt.pb.2} to a future work.
\begin{theo}[Convergence to the steady state] \label{theo.lt}
Let $(n^{eq}, V^{eq})\in H^1(\Omega)^4\times H^1(\Omega)$ satisfy eqs.~\eqref{lt.hp}--\eqref{lt.hp.3}. A constant $K>0$ depending on 
$D$, $p$, $\Omega$ exists such that, if:
$$(v^{sat})^2 + L^2\lambda_D^{-4}\|n^{eq}\|_{L^\infty(\Omega)}^2 < K , $$ 
where $L$ is the constant in eq.~\eqref{hp.mu}, then:
\begin{equation}
 \|n(t) - n^{eq}\|_{L^2(\Omega)} \leq k_1 e^{-k_2 t}\qquad t>0 ,\label{lt.decay}
\end{equation}
for suitable constants $k_1$, $k_2>0$.
\end{theo}
The idea for the proof is to derive an entropy inequality for the relative entropy:
$$ S[n,n^{eq}] = \frac{1}{2}\int_\Omega |n-n^{eq}|^2 , $$
which reads:
\begin{align} 
 &\frac{d}{dt}\int_\Omega S[n,n^{eq}] + c_0\int_\Omega |\nabla(n-n^{eq})|^2\leq 
 c( (v^{sat})^2 + L^2\lambda_D^{-4}\|n^{eq}\|_{L^\infty(\Omega)}^2 )S[n,n^{eq}] . \label{intro.lt}
\end{align}
Poincar\'e Lemma implies that the entropy $S[n,n^{eq}]$ can be controlled by the entropy dissipation $\int_\Omega |\nabla(n-n^{eq})|^2$.
So, under a suitable smallness assumption on the constant $(v^{sat})^2 + L^2\lambda_D^{-4}\|n^{eq}\|_{L^\infty(\Omega)}^2$, the entropy
inequality \eqref{intro.lt} and Gronwall lemma imply eq.~\eqref{lt.decay}.

To prove the last Theorem, concerning the regularity of the solutions to eqs.~\eqref{eq.n}--\eqref{eq.ic}, we need the following result
(an immediate corollary of Theorem 1.1 in \cite{Amann}), which is a generalization of the well-known Aubin lemma.
\begin{prop}\label{prop.amann} 
 Let $E$, $E_0$, $E_1$ Banach spaces such that the embedding $E_1\hookrightarrow E_0$ is compact and the embeddings
 $E_1\hookrightarrow E$, $E\hookrightarrow E_0$ are continuous. Moreover let us assume that constants $c>0$, $\theta\in (0,1)$ exist such that:
 \begin{align}
  \|u\|_E\leq c\|u\|_{E_0}^{1-\theta}\|u\|_{E_1}^\theta\qquad\forall u\in E_1 .\label{interp}
 \end{align}
 Then, for all $p\in [1,\infty]$, the embedding $L^p(0,T; E_1)\cap W^{1,p}(0,T; E_0)\hookrightarrow C([0,T]; E)$ is compact.
\end{prop}
\begin{theo}[Improved regularity]\label{theo.ireg} Let $\Omega\subset\R^2$,
 $\lambda_D>0$, $\gamma>0$, $T>0$. Moreover let us assume that $\vec m$ satisfies eq.~\eqref{reg.m.2} and
\begin{align}
 & n^{I}_j\in H^1(\Omega),\quad n_j^D,\, V^D\in H^{2}(\Omega),\quad C\in L^\infty(\Omega)\qquad (0\leq j\leq 3),\label{ireg.hp.data}\\
 & D, p\in W^{1,\infty}(\Omega),\quad \inf_\Omega D>0,\quad \sup_\Omega |p| < 1.\label{ireg.hp.parms}
\end{align}
Then the unique solution $(n,V)$ to problem \eqref{eq.n}--\eqref{eq.ic} satisfies:
\begin{align}
 & n_j\in L^q(0,T; W^{2,q}(\Omega))\cap W^{1,q}(0,T; L^q(\Omega))\cap C([0,T]; W^{1,q}(\Omega)) 
 \qquad (0\leq j\leq 3) ,\label{reg.n.2}\\
 & V \in L^q(0,T; W^{3,q}(\Omega))\cap W^{1,q}(0,T; W^{1,q}(\Omega))\cap C([0,T]; W^{2,q}(\Omega)), \label{reg.V.2}
\end{align}
for all $q\in [1,2)$.
\end{theo}
The idea of the proof is that the regularity assumptions \eqref{reg.m.2} on $\vec m$, stronger than the hypothesis \eqref{reg.m} 
made for the existence analysis, allow us 
to derive and exploit a set of equations for the variables $n_+$, $n_-$, $\vec n_\perp$ defined in eq.~\eqref{newvar}.
Such equations are a generalization of eqs. (9), (10), (13) in \cite{JueNegShp} and have the form:
\begin{align}
 & \pa_t n_\pm - \pa_i\left(\frac{D}{1\pm p}(\pa_i n_\pm - v_i n_\pm)\right) = f_\pm, \quad
 \pa_t\vec n_\perp - \pa_i\left(\frac{D}{\eta}(\pa_i\vec n_\perp - v_i\vec n_\perp)\right) 
 + \frac{\vec n_\perp}{\tau} = \vec f_\perp ,\label{intro.f}
\end{align}
with $f_\pm$, $\vec f_\perp$ suitable quantities depending from $n$, $\nabla n$, $\vec m$, $\pa_t\vec m$, $\nabla\vec m$, $\Delta\vec m$, $v$.
The big advantage in eqs.~\eqref{intro.f} consists in
the lack of cross diffusion terms, which makes possible to exploit the above equations in order to derive improved regularity results for the
solutions of pb.~\eqref{eq.n}--\eqref{eq.ic}.
The proof of the Theorem consists in three parts. In the first part eq.~\eqref{intro.f} is derived.
In the second part eqs.~\eqref{reg.n.2}, \eqref{reg.V.2} are shown for $p=4/3$.
Hypothesis $\Omega\subset\R^2$ allows to employ the following Gagliardo-Nirenberg inequality \cite{Brezis}:
\begin{align}
 \|u\|_{L^4}\leq c \|u\|_{L^2}^{1/2}\|u\|_{H^1}^{1/2}\qquad u\in H^1(\Omega),\label{intro.GN}
\end{align}
which, together with the already known regularity results \eqref{reg.n} on $n$, implies that $n \in L^4(\OmT)$. 
From this fact and the assumptions \eqref{reg.m.2} on $\vec m$, we are able to deduce
that the source terms $f_\pm$, $\vec f_\perp$ in eq.~\eqref{intro.f} belong to $L^{4/3}(\OmT)$. This property, thanks to suitable regularity results
for parabolic problems \cite{Brezis, Fri, LadSolUra}, yields eqs.~\eqref{reg.n.2}, \eqref{reg.V.2} for $p=4/3$.
In the third part of the proof, eqs.~\eqref{reg.n.2}, \eqref{reg.V.2} are shown to hold for an increasing sequence of exponents 
$(p_k)_{k\in\N}\subset (1,2)$.
The partial result obtained in the previous part of the proof provides the first step in the derivation of such a sequence ($p_0=4/3$), which is
built by iteratively applying Proposition \ref{prop.amann} with $E_0 = L^{p_k}(\Omega)$, $E = W^{1,p_k}_0(\Omega)$, 
$E_1 = W^{2,p_k}\cap W^{1,p_k}_0(\Omega)$. In particular, assumption \eqref{interp} is a consequence of the following Gagliardo-Nirenberg inequality
\cite{Brezis}:
\begin{equation}
 \|\nabla u\|_{L^{p}}\leq c \|u\|_{L^{p}}^{1/2}\|u\|_{W^{2,p}}^{1/2}\qquad u\in W^{2,p}(\Omega),\label{intro.GN.2}
\end{equation}
valid for $p\in [1,\infty]$. The fact that $p_k\to 2$ as $k\to\infty$ implies that eqs.~\eqref{reg.n.2}, \eqref{reg.V.2} hold for all $p<2$. 

The paper is organized as follows.
We prove Theorem \ref{theo.ex} in Section \ref{sec.ex}. Theorem \ref{theo.lt} is proved in Section \ref{sec.lt}. Finally,
Section \ref{sec.ireg} is devoted to the proof of Theorem \ref{theo.ireg}.

\section{Existence and uniqueness of weak solutions}\label{sec.ex}
In this section we prove Theorem \ref{theo.ex}. The proof is divided into four steps.

{\em Step 1: time semi-discretization.} 
Let $h>0$ be the time step, $t_k = kh$, $n^{(k)} = n(\cdot,t_k)$, $V^{(k)} = V(\cdot,t_k)$ for $0\leq k\leq T/h$.
We consider the following semi-discretized problem: 
\begin{align}
 & n_i^{(k)} - h\Div (a_{ij}^{(k)}\nabla n_j^{(k)}) 
 = n_i^{(k-1)} - h\Div (a_{ij}^{(k)}n_j^{(k)}v^{(k)})
 + h b_{ij}^{(k)}n_j^{(k)}\quad (0\leq i\leq 3)\quad\textrm{in }\Omega,\label{eq.n.d} \\
 & v^{(k)} = -\mu(|\nabla V^{(k)}|)\nabla V^{(k)}\quad\textrm{in }\Omega,\label{eq.v.d}\\
 & -\lambda_D^2\Delta V^{(k)} = 2 n_0^{(k)} - C\qquad\textrm{in }\Omega,\label{eq.poi.d}\\
 & n^{(k)} = n^D,\quad V^{(k)}=V^D\qquad\textrm{on }\pa\Omega,\label{eq.bc.d}
\end{align}
for $k\geq 1$.

{\em Step 2: fixed point.} We solve eqs.~\eqref{eq.n.d}--\eqref{eq.bc.d} by applying Leray-Schauder's fixed point theorem \cite{Zei90}.
Let us define the operator 
$$\fun{F}{L^2(\Omega)\times [0,1]}{L^2(\Omega)}{(n,\sigma)}{u} $$
where $u\in H^1(\Omega)$ satisfies:
\begin{align}
 & u_i - h\Div (a_{ij}^{(k)}\nabla u_j) 
 = n_i^{(k-1)} - \sigma h\Div (a_{ij}^{(k)} n_j v)
 + \sigma h b_{ij}^{(k)}n_j\qquad (0\leq i\leq 3)\quad\textrm{in }\Omega,\label{eq.n.fp}\\
 & v = -\mu(|\nabla V|)\nabla V\quad\textrm{in }\Omega,\label{eq.v.fp}\\
 & -\lambda_D^2\Delta V = 2 n_0 - C\qquad\textrm{in }\Omega,\label{eq.poi.fp}\\
 & u = n^D,\quad V=V^D\qquad\textrm{on }\pa\Omega.\label{eq.bc.fp}
\end{align}
The existence of a unique solution $u\in H^1(\Omega)$ to eqs.~\eqref{eq.n.fp}--\eqref{eq.bc.fp} is an easy application of Lax-Milgram lemma \cite{Brezis}.
The operator $F$ is compact due to the compact Sobolev embedding $H^1(\Omega)\hookrightarrow L^2(\Omega)$. 
The continuity of $F$ can be proved with a standard argument (see e.g.~\cite{Jue94}). Moreover $F(\cdot,0)$ is constant.
Finally, let $n\in L^2(\Omega)$, $\sigma\in [0,1]$ such that $F(n,\sigma)=n$. Let us consider the quadratic entropy functional:
\begin{equation}\label{Phi}
 S[n] = \frac{1}{2}\int_\Omega |n-n^D|^2 .
\end{equation}
From eq.~\eqref{eq.n.fp} and the fact that $A$ is symmetric and positive definite it follows:
\begin{align}
 & h^{-1}(S[n]-S[n^{(k-1)}]) \leq \int_\Omega (n_j-n_j^D)h^{-1}(n_j-n_j^{(k-1)})\nonumber\\
 &\quad = -\int_\Omega (a_{ij}^{(k)}\nabla n_j - \sigma a_{ij}^{(k)}n_j v)\cdot\nabla (n_i-n_i^D)
 + \int_\Omega \sigma b_{ij}^{(k)}n_j (n_i-n_i^D)\nonumber\\
 &\quad\leq -\frac{1}{2}\int_\Omega a_{ij}^{(k)}\nabla n_i\cdot\nabla n_j + \frac{1}{2}\int_\Omega a_{ij}^{(k)}\nabla n_i^D\cdot\nabla n_j^D 
 + c\int_\Omega |n|(|\nabla(n-n^D)| + |n-n^D|),\nonumber
\end{align}
so by exploiting the positivity and boundedness of $A$ and the assumptions \eqref{ex.hp.1} on the data and by applying Young inequality we obtain:
\begin{align}
 & h^{-1}(S[n]-S[n^{(k-1)}]) + c_0\int_\Omega |\nabla n|^2 \leq c(S[n]+1),\label{dS.fp}
\end{align}
for some positive constants $c_0$, $c$, independent on $\sigma$, $h$.\\
Eq.~\eqref{dS.fp} yields $\sigma-$uniform bounds for $n$, $\nabla n$ in $L^2(\Omega)$.
So from Leray-Schauder's theorem (see e.g.~\cite{Zei90}) we obtain the existence of a fixed point 
$n^{(k)}\in H^1(\Omega)$ for $F(\cdot,1)$, which means, a solution of
\eqref{eq.n.d}--\eqref{eq.bc.d}.

{\em Step 3: limit $h\to 0$.} Eq.~\eqref{dS.fp} holds with $n=n^{(k)}$:
\begin{align}
 & h^{-1}(S[n^{(k)}]-S[n^{(k-1)}]) + c_0\int_\Omega |\nabla n^{(k)}|^2 \leq c(S[n^{(k)}]+1).\label{dS.discr}
\end{align}
Let us define the piecewise constant functions: $n^h(\cdot,t) = n^{(k)}$, $V^h(\cdot,t) = V^{(k)}$ for $t\in ((k-1)h,kh]$.  
From eq.~\eqref{dS.discr} we deduce:
\begin{align}
 & S[n^h(t)] + c_0\int_0^t\int_\Omega |\nabla n^{h}|^2 \leq c_1\int_0^t S[n^h(s)]\,ds + c(T),
 \qquad t\in [0,T].\label{dS}
\end{align}
By applying Gronwall's lemma, from eq.~\eqref{dS} we get:
\begin{equation}
 \|n^h\|_{L^2(0,T; H^1(\Omega))} + \|n^h\|_{L^\infty(0,T; L^2(\Omega))}\leq c .\label{est.nh}
\end{equation}
If we define the discrete time derivative $D_h w$ of an arbitrary function $u\,:\,\OmT\to\R$ as: 
$$D_h w(x,t) = (w(x,t)-w(x,t-h))/h \qquad (x,t)\in\Omega\times (h,T), $$
from eqs.~\eqref{eq.n.d}, \eqref{est.nh} it follows easily:
\begin{equation}
 \|D_h n^h\|_{L^2(h,T; H^{-1}(\Omega))}\leq c .\label{est.dnh}
\end{equation}
By exploiting eqs.~\eqref{est.nh}, \eqref{est.dnh} and the compact embedding $H^1(\Omega)\hookrightarrow L^p(\Omega)$ for all $p\in [1,6)$,
from Theorem 1 in \cite{DJ11} we argue that, up to subsequences, 
$n^h\to n$ in $L^2(0,T; L^p(\Omega))$ as $h\to 0$, for all $p\in [1,6)$. Since eq.~\eqref{eq.poi.d} holds, this implies
$V^h\to V$ in $L^2(0,T; W^{2,p}(\Omega))$ as $h\to 0$; in particular, $\nabla V^h\to\nabla V$ a.e. in $\OmT$ as $h\to 0$, so from the dominated
convergence theorem $v^k\to v = -\mu(|\nabla V|)\nabla V$ in $L^q(\OmT)$ for all $q\in [1,\infty)$.
From the pointwise convergence of $n^h$ and Fatou's lemma we deduce that $n\in L^\infty(0,T; L^2(\Omega))$.
By exploiting these properties it is straightforward to take the limit $h\to 0$ in 
eqs.~\eqref{eq.n.d}--\eqref{eq.bc.d}, proving that $(n, V)$ satisfies eqs.~\eqref{eq.n}--\eqref{eq.ic}, \eqref{reg.n}, \eqref{reg.V}.

{\em Step 4: uniqueness.} 
Let $(n, V)$, $(\bn, \bV)$ solutions to eqs.~\eqref{eq.n}--\eqref{eq.ic} satisfying eqs.~\eqref{reg.n}, \eqref{reg.V}.
Moreover let $v = -\mu(|\nabla V|)\nabla V$, $\overline{v}=-\mu(|\nabla \bV|)\nabla \bV$.
Then:
\begin{align}
 &\pa_t (n_i-\bn_i) = \Div (a_{ij}(\nabla (n_j-\bn_j) - (n_j-\bn_j)v
 - \bn_j (v-\overline v)) )\nonumber\\
 & \quad + b_{ij}(n_j-\bn_j)\quad (0\leq i\leq 3)\qquad\textrm{in }\Omega\times (0,T),\label{eq.n.uni}\\
 & -\lambda_D^2\Delta (V-\bV) = n_0 - \bn_0\qquad\textrm{in }\OmT,\label{eq.poi.uni}\\
 & n-\bn = 0,\quad V-\bV=0\qquad\textrm{on }\pa\OmT,\label{eq.bc.uni}\\
 & n(\cdot,0)-\bn(\cdot,0) = 0\qquad\textrm{on }\Omega.\label{eq.ic.uni}
\end{align}
Because of eq.~\eqref{hp.mu} it holds that $|v-\overline{v}|\leq L |\nabla(V-\bV)|$.
So, if we use $n - \bn$ as a test function in the weak formulation of eq.~\eqref{eq.n.uni} we find:
\begin{align}
 & \frac{d}{dt}\int_\Omega\frac{1}{2}|n-\bn|^2 + c_0 \int_\Omega |\nabla(n-\bn)|^2 \label{uni.1}\\
 &\quad \leq c \int_\Omega |n-\bn| |\nabla(n-\bn)| + c \int_\Omega |\bn| |\nabla (V-\bV)| |\nabla(n-\bn)|\nonumber\\
 &\quad \leq\ep\int_\Omega |\nabla(n-\bn)|^2 
 + \frac{c}{\ep}\int_\Omega\left(|n-\bn|^2 + |\bn|^2 |\nabla (V-\bV)|^2 \right)\nonumber\\
 &\quad \leq\ep\int_\Omega |\nabla(n-\bn)|^2 + \frac{c}{\ep}(\|n-\bn\|_{L^2(\Omega)}^2
 + \|\bn\|_{L^4(\Omega)}^2 \|\nabla (V-\bV)\|_{L^4(\Omega)}^2 ) . \nonumber
\end{align}
From the Sobolev embedding $H^1(\Omega)\subset L^4(\Omega)$ and eq.~\eqref{eq.poi.uni} it follows:
$$ 
\|\bn\|_{L^4(\Omega)}^2 \|\nabla (V-\bV)\|_{L^4(\Omega)}^2 
\leq c\|\bn\|_{H^1(\Omega)}^2 \|V-\bV\|_{H^2(\Omega)}^2
\leq c\|\bn\|_{H^1(\Omega)}^2 \|n-\bn\|_{L^2(\Omega)}^2 ,
$$
so eq.~\eqref{uni.1} can be rewritten as:
\begin{align}
 & \frac{d}{dt}\int_\Omega\frac{1}{2}|n-\bn|^2 + c_0 \int_\Omega |\nabla(n-\bn)|^2 
 \leq\ep\int_\Omega |\nabla(n-\bn)|^2 
 + \frac{c}{\ep}(1+\|\bn\|_{H^1(\Omega)}^2)\|n-\bn\|_{L^2(\Omega)}^2 . \nonumber
\end{align}
If we choose $\ep>0$ small enough we obtain:
\begin{equation}
 \frac{d}{dt}\|n-\bn\|_{L^2(\Omega)}^2 \leq c(1+\|\bn\|_{H^1(\Omega)}^2)\|n-\bn\|_{L^2(\Omega)}^2\qquad t>0.\label{uni.4}
\end{equation}
Eq.~\eqref{reg.n} implies that $\|\bn\|_{H^1(\Omega)}^2\in L^1(0,T)$, thus
from Gronwall's inequality we deduce that $n=\bn$ in $\OmT$. 
From eqs.~\eqref{eq.poi.uni}, \eqref{eq.bc.uni} we also get $V=\bV$ in $\OmT$. This finishes the proof.

\section{Long-time behaviour of solutions}\label{sec.lt}

In this section we prove Theorem \ref{theo.lt}.

Let us employ $n-n^{eq}$ as test function in the weak formulation of eqs.~\eqref{eq.n}--\eqref{eq.ic}. We obtain:
\begin{align} 
 &\frac{d}{dt}\int_\Omega\frac{|n-n^{eq}|^2}{2} + \int_\Omega\pa_i(n-n^{eq})\cdot A\pa_i n 
  - \int_\Omega v_i\pa_i(n-n^{eq})\cdot A n - \int_\Omega (n-n^{eq})\cdot B n = 0,\nonumber
\end{align}
which, by applying eq.~\eqref{lt.hp}, can be rewritten as:
\begin{align} 
 &\frac{d}{dt}\int_\Omega\frac{|n-n^{eq}|^2}{2} + \int_\Omega\pa_i(n-n^{eq})\cdot A\pa_i (n-n^{eq}) \nonumber\\
 &\quad - \int_\Omega \pa_i(n-n^{eq})\cdot A ( (n-n^{eq}) v_i + n^{eq} (v_i - v_i^{eq}) ) - \int_\Omega (n-n^{eq})\cdot B n = 0.\label{lt.2}
\end{align}
Since $A$ is strictly positive definite we deduce:
\begin{align}
 &\int_\Omega\pa_i(n-n^{eq})\cdot A\pa_i (n-n^{eq})\geq c_0\int_\Omega |\nabla(n-n^{eq})|^2 ,
 \label{lt.3}
\end{align}
for some $c_0 > 0$. Moreover, since $|v|\leq v^{sat}$:
\begin{align}
 &\int_\Omega \pa_i(n-n^{eq})\cdot A (n - n^{eq})v_i
 \leq c\ep\int_\Omega |\nabla(n-n^{eq})|^2 + \frac{c}{\ep}(v^{sat})^2\int_\Omega |n-n^{eq}|^2 . \label{lt.4}
\end{align}
From eqs.~\eqref{eq.poi.bis}, \eqref{eq.bc}, \eqref{lt.poi}, \eqref{lt.hp.2} it follows:
\begin{align}
 & -\lambda_D^2\Delta(V-V^{eq}) = n_0 - n_0^{eq}\qquad\textrm{in }\OmT,\label{lt.dpoi}\\
 & V-V^{eq} = 0\qquad\textrm{on }\pa\OmT.\label{lt.dpoi.bc}
\end{align}
Thus eq.~\eqref{hp.mu} and standard estimates for linear elliptic equations applied to eqs.~\eqref{lt.dpoi}, \eqref{lt.dpoi.bc}
imply:
\begin{align}
 & \int_\Omega \pa_i(n-n^{eq})\cdot A n^{eq} (v_i - v_i^{eq})\nonumber\\
 & \quad\leq c\ep\int_\Omega |\nabla(n-n^{eq})|^2 + \frac{c}{\ep} L^2\|n^{eq}\|_{L^\infty(\Omega)}^2\int_\Omega |\nabla(V-V^{eq})|^2 \nonumber\\
 & \quad\leq c\ep\int_\Omega |\nabla(n-n^{eq})|^2 + \frac{c}{\ep} L^2\lambda_D^{-4}\|n^{eq}\|_{L^\infty(\Omega)}^2\int_\Omega |n-n^{eq}|^2 .\label{lt.5}
\end{align}
Finally, from the definitions of $B$ and $n^{eq}$ it follows immediately:
\begin{align}
 & -\int_\Omega (n-n^{eq})\cdot B n = \int_\Omega\frac{|\vec n|^2}{\tau}\geq 0 . \label{lt.6}
\end{align}
By collecting eqs.~\eqref{lt.2}--\eqref{lt.4}, \eqref{lt.5}, \eqref{lt.6} and choosing $\ep>0$ small enough we find:
\begin{align} 
 &\frac{d}{dt}\int_\Omega |n-n^{eq}|^2 + c_0\int_\Omega |\nabla(n-n^{eq})|^2\leq 
 c( (v^{sat})^2 + L^2\lambda_D^{-4}\|n^{eq}\|_{L^\infty(\Omega)}^2 )\int_\Omega |n_0-n_0^{eq}|^2 . \label{lt.7}
\end{align}
It is straightforward fo deduce that the constants $c_0$, $c$ in eq.~\eqref{lt.7} depend only on $D$ and $p$.
Poincar\'e's Lemma implies that $\int_\Omega |\nabla(n-n^{eq})|^2\geq c_1\int_\Omega |n-n^{eq}|^2$, for a suitable constant $c_1>0$ depending only 
on $\Omega$. So, if $(v^{sat})^2 + L^2\lambda_D^{-4}\|n^{eq}\|_{L^\infty(\Omega)}^2 < K\equiv c_0 c_1/c$
then eq.~\eqref{lt.decay} follows from Gronwall's Lemma applied to eq.~\eqref{lt.7}. This finishes the proof.

\section{Improved regularity}\label{sec.ireg}
In this section we prove Thorem \ref{theo.ireg}. The proof is divided into three steps.

{\em Step 1: diagonalization of the system.} 
By exploiting eqs.~\eqref{eq.n}--\eqref{eq.ic}, we derive a set of equations, without cross diffusion terms, 
for the variables $n_+$, $n_-$, $\vec n_\perp$ defined in eq.~\eqref{newvar}.
Let us first multiply eq.~\eqref{eq.nk} times $m_k$ and sum over $k=1,2,3$. By exploiting the fact that $\pa_i\vec m\cdot\vec m = 0$ for $i=1,2$
(thanks to eq.~\eqref{LLG.eq}) we get:
\begin{align}
 & \pa_t(\vec n\cdot\vec m) - \vec n\cdot\pa_t\vec m - \pa_i\left(\frac{D}{\eta^2}(-p J^i_0 + m_s J^i_s)\right)
 + \frac{D}{\eta}J_s^i\pa_i m_s + \frac{\vec n\cdot\vec m}{\tau} = 0.\label{ireg.1}
\end{align}
Now, if we compute the sum and the difference between eqs.~\eqref{eq.n0}, \eqref{ireg.1}, we obtain equations for $\pa_t n_+$, $\pa_t n_-$,
respectively:
\begin{align}
 & \pa_t n_\pm \mp\vec n\cdot\pa_t\vec m - \pa_i\left(\frac{D}{1\pm p}(J_0^i\pm m_s J_s^i)\right)
 \pm \frac{D}{\eta}J_s^i\pa_i m_s \pm \frac{\vec n\cdot\vec m}{\tau} = 0.\label{ireg.2}
\end{align}
By exploiting the simple relations:
\begin{align}
 J_0^i\pm m_s J_s^i = \pa_i n_\pm - v_i n_\pm \mp\vec n\cdot\pa_i\vec m , \quad
 J_s^i\pa_i m_s = \pa_i \vec m\cdot\pa_i \vec n - v_i\vec n\cdot\pa_i\vec m ,\nonumber
\end{align}
we can rewrite eq.~\eqref{ireg.2} as:
\begin{align}
 & \pa_t n_\pm - \pa_i\left(\frac{D}{1\pm p}(\pa_i n_\pm - v_i n_\pm)\right) = 
 \pm\vec n\cdot\pa_t\vec m \mp \pa_i\left(\frac{D}{1\pm p}\vec n\cdot\pa_i\vec m\right)\nonumber\\
 &\quad \mp\frac{D}{\eta}(\pa_i \vec m\cdot\pa_i \vec n - v_i\vec n\cdot\pa_i\vec m) \mp \frac{\vec n\cdot\vec m}{\tau}. \label{eq.npm}
\end{align}
Now let us multiply eq.~\eqref{ireg.1} times $m_k$:
\begin{align}
 & m_k\pa_t(\vec n\cdot\vec m) - m_k\vec n\cdot\pa_t\vec m - \pa_i\left(\frac{D}{\eta^2}m_k(-p J^i_0 + m_s J^i_s)\right)\nonumber\\
 &\quad + \frac{D}{\eta^2}(-p J^i_0 + m_s J^i_s)\pa_i m_k
 + \frac{D}{\eta}J_s^i m_k\pa_i m_s + m_k\frac{\vec n\cdot\vec m}{\tau} = 0 \qquad(k=1,2,3).\label{ireg.3}
\end{align}
Since: 
$$\pa_t(n_\perp)_k = \pa_t n_k - m_k\pa_t(\vec n\cdot\vec m) - (\vec n\cdot\vec m)\pa_t m_k , $$
by taking the difference between eqs.~\eqref{eq.nk}, \eqref{ireg.3} we get:
\begin{align}
 & \pa_t(n_\perp)_k + (\vec n\cdot\vec m)\pa_t m_k + m_k\vec n\cdot\pa_t\vec m - \pa_i\left(\frac{D}{\eta}(\delta_{ks}-m_k m_s)J^i_s\right)
 -2\gamma\ep_{ijk}n_i m_j\nonumber\\
 &\quad + \frac{D}{\eta^2}(p J^i_0 - m_s J_s^i)\pa_i m_k - \frac{D}{\eta}J^i_s m_k\pa_i m_s + \frac{(n_\perp)_k}{\tau} = 0
 \qquad(k=1,2,3).  \label{ireg.4}
\end{align}
Elementary computations imply that:
\begin{align}
 & (\delta_{ks}-m_k m_s)J^i_s = \pa_i(n_\perp)_k - v_i(n_\perp)_k + n_s\pa_i(m_k m_s),\nonumber\\
 & (\vec n\cdot\vec m)\pa_t m_k + m_k\vec n\cdot\pa_t\vec m = n_s\pa_t(m_k m_s), \nonumber
\end{align}
so from eq.~\eqref{ireg.4} we conclude:
\begin{align}
 & \pa_t(n_\perp)_k - \pa_i\left(\frac{D}{\eta}(\pa_i(n_\perp)_k - v_i(n_\perp)_k)\right) + \frac{(n_\perp)_k}{\tau} = 
 \pa_i\left(\frac{D}{\eta}n_s\pa_i(m_k m_s)\right)\nonumber\\
 &\quad -\frac{D}{\eta^2}( p\delta_{k\ell}(\pa_i n_0 - v_i n_0) 
   -m_\alpha(\delta_{\alpha\beta}\delta_{k\ell} + \eta\delta_{\alpha k}\delta_{\beta\ell})(\pa_i n_\beta - v_i n_\beta)
 )\pa_i m_\ell \nonumber\\
 &\quad -n_s\pa_t(m_k m_s) + 2\gamma\ep_{ijk}n_i m_j \qquad (k=1,2,3).\label{eq.nu}
\end{align}

{\em Step 2: the thesis holds for $p=4/3$.} 
We look now for estimates of the right-hand sides of eqs.~\eqref{eq.npm}, \eqref{eq.nu}.
In order to do this, let us observe that, since $\Omega\subset\R^2$, we can apply the following Gagliardo-Nirenberg inequality \cite{Brezis},
valid for $1\leq q\leq p < \infty$:
\begin{align}
 & \|n_s(t)\|_{L^p(\Omega)}\leq c\|n_s(t)\|_{L^q(\Omega)}^{q/p}\|n_s(t)\|_{H^1(\Omega)}^{1-q/p}\quad (0\leq s\leq 3)\quad\textrm{a.e. }t\in (0,T),
 \label{ireg.GN}
\end{align}
where $c>0$ is a suitable constant, depending only on $\Omega$. Since eq.~\eqref{reg.n} holds, then eq.~\eqref{ireg.GN} with 
$p=4$, $q=2$ implies:
\begin{align}
 & \|n_s(t)\|_{L^4(\Omega)}^4 \leq \tilde{c}_s\|n_s(t)\|_{H^1(\Omega)}^{2}\quad (0\leq s\leq 3)\quad\textrm{a.e. }t\in (0,T),
 \label{ireg.GN.2}
\end{align}
where the constant $\tilde c_s>0$ is proportional to $\|n_s\|_{L^\infty(0,T; L^2(\Omega))}^2<\infty$. 
By integrating eq.~\eqref{ireg.GN.2} in time and exploiting eq.~\eqref{reg.n} we deduce that $n_s\in L^4(\OmT)$, for $0\leq s\leq 3$.
We are going to exploit this fact and eq.~\eqref{reg.m.2} to prove that the right-hand sides of eqs.~\eqref{eq.npm}, \eqref{eq.nu} belong to $L^{4/3}(\OmT)$. 
Each term appearing in the right-hand sides of eqs.~\eqref{eq.npm}, \eqref{eq.nu} can be bounded (in modulus) by the function:
\begin{align}
 & f \equiv c|n|(|\pa_t\vec m| + |\Delta \vec m| + |\nabla\vec m|^2 + 1) + c|\nabla n||\nabla\vec m| ,
 \label{ireg.bad}
\end{align}
with $c>0$ a suitable constant. From eq.~\eqref{reg.m.2} it follows:
\begin{align}
 & |\pa_t\vec m| + |\Delta \vec m| + |\nabla\vec m|^2 + 1 \in L^\infty(0,T; L^2(\Omega)),
 \quad |\nabla\vec m|\in L^\infty(0,T; L^p(\Omega))\quad\forall p<\infty,
 \nonumber
\end{align}
which, together with eq.~\eqref{reg.n} and the fact that $n_s\in L^4(\OmT)$ for $0\leq s\leq 3$, implies $f\in L^{4/3}(\OmT)$, so 
also the right-hand sides of eqs.~\eqref{eq.npm}, \eqref{eq.nu} belong to $L^{4/3}(\OmT)$. 
From this fact and the hypothesis \eqref{ireg.hp.data} on the data, by applying standard regularity results for parabolic equations 
(see e.g.~\cite{Brezis, Fri, LadSolUra}) it follows that:
\begin{align}
 & n_\pm, (n_\perp)_k \in L^{4/3}(0,T; W^{2,4/3}(\Omega))\cap W^{1,4/3}(0,T; L^{4/3}(\Omega))\qquad (k=1,2,3).
 \label{ireg.npmnu}
\end{align}
By exploiting eqs.~\eqref{reg.m.2}, \eqref{ireg.npmnu} and the inverse relations
$n_0 = (n_+ + n_-)/2$, $\vec n = (n_+ - n_-)\vec m/2 + \vec n_\perp$, we get:
\begin{align}
 & n_s,\, \pa_{t}n_s,\, \pa_{x_i}n_s,\, \pa_{x_i x_j}^2 n_s \in L^{4/3}(\OmT)
 \qquad 0\leq s\leq 3,\,\, 1\leq i,j\leq 2 , \label{reg.n.p0}
\end{align}
namely eq.~\eqref{reg.n.2} for $p=4/3$. Eq.~\eqref{reg.V.2} for $p=4/3$ is a consequence of eqs.~\eqref{eq.poi}, \eqref{reg.n.p0}.

{\em Step 3: the thesis holds for $p<2$.} 
Now we prove the following:
\begin{claim}
An increasing sequence $(p_k)_{k\in\N}\subset (1,2)$ exists such that:
\begin{align}
 & n_s,\, \pa_{t}n_s,\, \pa_{x_i}n_s,\, \pa_{x_i x_j}^2 n_s \in L^{p_k}(\OmT)
 \qquad 0\leq s\leq 3,\,\, 1\leq i,j\leq 2,\label{reg.n.pk}
\end{align}
for all $k\geq 0$.
\end{claim}
{\em Proof.}
The proof is by induction on $k$. Let $p_0 \equiv 4/3$, so the thesis is true for $k=0$ as a consequence of Step 1.
Let us assume eq.~\eqref{reg.n.pk} is true for some $k\geq 0$. Let us define the following reflexive Banach spaces:
\begin{align}
 E_1 = W^{2,p_k}\cap W^{1,p_k}_0(\Omega),\qquad
 E = W^{1,p_k}_0(\Omega),\qquad
 E_0 = L^{p_k}(\Omega).\nonumber
\end{align}
From standard Sobolev embedding theorems it follows that the embeddings $E_1\hookrightarrow E$, $E\hookrightarrow E_0$ 
are compact. Moreover Poincar\'e inequality and eq.~\eqref{intro.GN.2} with $p=p^k$ imply:
\begin{align}
 \|u\|_E \leq c\|\nabla u\|_{L^{p_k}} \leq c\|u\|_{E_0}^{1/2}\|u\|_{E_1}^{1/2}\qquad\forall u\in E_1 .
\end{align}
With this notation, eq.~\eqref{reg.n.pk} means that $n_s-n_s^D\in L^{p_k}(0,T; E_1)$ and $\pa_t(n_s-n_s^D)\in L^{p_k}(0,T; E_0)$ for $0\leq s\leq 3$.
Thus, by applying Prop.~\ref{prop.amann}, we deduce that $n_s-n_s^D \in C([0,T]; E)$, which implies,
again by means of standard Sobolev embedding theorems, that: 
\begin{align}
& n_s\in C([0,T]; W^{1,p_k}(\Omega))\hookrightarrow 
L^\infty(0,T; L^{p_k^*}(\Omega))\quad (0\leq s\leq 3), \qquad p_k^* \equiv \frac{2 p_k}{2-p_k}  . \label{reg.n.inf}
\end{align}
From eqs.~\eqref{reg.n}, \eqref{ireg.GN} with $q=p_k^*$, $p = p_k^* + 2$ and \eqref{reg.n.inf} it follows that 
$n_s\in L^{p_k^* + 2}(\OmT)$. By arguing like in the derivation of eq.~\eqref{reg.n.p0} in step 2 we deduce that:
\begin{align}
 & n_s,\, \pa_{t}n_s,\, \pa_{x_i}n_s,\, \pa_{x_i x_j}^2 n_s \in L^{p_{k+1}}(\OmT)
 \qquad 0\leq s\leq 3,\,\, 1\leq i,j\leq 2,\label{reg.n.pk1}
\end{align}
with $p_{k+1} = \left( \frac{1}{2} + \frac{1}{p_k^* + 2} \right)^{-1}$, which means:
\begin{align}
 & p_{k+1} = \frac{4}{4-p_k}\qquad \forall k\geq 0,\qquad p_0 = \frac{4}{3} . \label{pk}
\end{align}
By exploiting eq.~\eqref{pk} it is easy to prove that $1 < p_k < p_{k+1} < 2$ for all $k\geq 0$. Thus the claim has been proved.
{\flushright $\Box$ \flushleft}
From eqs.~\eqref{eq.poi}, \eqref{reg.n.pk}, \eqref{reg.n.inf} it follows also:
\begin{align}
 & V \in L^{p_k}(0,T; W^{3,{p_k}}(\Omega))\cap W^{1,{p_k}}(0,T; W^{1,{p_k}}(\Omega))\cap C([0,T]; W^{2,{p_k}}(\Omega))
 \qquad (k\geq 0). \label{reg.V.pk}
\end{align}
The monotonicity of $p_k$ implies that $p_k$ has a finite limit $\tilde p$ which must satisfy $\tilde p = 4/( 4 - \tilde p )$, so $\tilde p = 2$. 
From this fact and eqs.~\eqref{reg.n.pk}, \eqref{reg.n.inf} and \eqref{reg.V.pk} we conclude that eqs.~\eqref{reg.n.2}, \eqref{reg.V.2} 
hold for all $q\in [1,2)$. This finishes the proof of Thorem \ref{theo.ireg}.

\section{Acknowledgements.} The author acknowledges economic support from the Austrian Science Fund (FWF) and thanks 
Univ.-Prof.~Dr.~Ansgar J\"ungel for the help received during the research work and the correction of the present paper.

\end{document}